\documentclass[12pt]{iopart}
\def\qed{\hfill $\square$}

\usepackage{iopams}  
\usepackage{graphicx,hyperref}
\usepackage{multirow,pifont,xcolor,amssymb} 
\usepackage{makecell}

\newtheorem{proposition}{Proposition}

\usepackage{xcolor}
\usepackage{placeins}

\begin{document}

\title[A dimension reduction procedure for the design of optimal lattice-spring systems]{A dimension reduction procedure for the design of lattice-spring systems with minimal fabrication cost and required multi-functional properties}

\author{Egor Makarenkov$^a$, Sakshi Malhotra$^b$, Yang Jiao$^c$}

\address{$^a$Allen High School, TX, USA}
\address{$^b$Mathematical Sciences,  University of Texas at Dallas, TX, USA}
\address{$^c$Material Sciences and Engineering Program, School for Engineering of Matter, Transport and Energy, Arizona State University, AZ, USA}
\ead{egordaporg@gmail.com, egor.makarenkov@student.allenisd.org, sakshi.malhotra@utdallas.edu, yang.jiao.2@asu.edu}
\vspace{10pt}
\begin{indented}
\item[]August 2017 (minor update March 2024)
\end{indented}

\begin{abstract}
We show that the problem of the design of the lattices of elastoplastic current conducting springs with optimal multi-functional properties leads to an analytically tractable problem.   Specifically, focusing on a lattice with a small number of springs, we use the technique of inequalities to reduce the number variables and to compute the minimal cost of lattice fabrication  explicitly. 
\end{abstract}

%

\section{Introduction} The problem of the design of optimal elastoplastic materials attracted lots of interest in the recent engineering literature \cite{topopt1,topopt2,topopt3,topopt4,topopt5}, but mathematical results mainly concern optimum existence theorems. The reason for that is the large dimension of the configuration space where explicit analytic computations are not feasible. Even though various dimension reduction techniques have been developed in optimization theory and in the field of the dynamics of elastoplastic materials, due to the complexity of the methods, the gap to analytic implementation of these methods is very big. 
The goal of this paper is to propose settings in which analytic computations towards optimization of elastoplastic material are doable and can be implemented by practitioners. 

\vskip0.2cm

\noindent Following Moreau \cite{ch5Moreau} we propose to view elastoplastic material as a lattice of connected elastoplastic springs. An elastoplastic spring has an ability to return to the original relaxed state if the elongation doesn't exceed certain critical value. The stress of an elastoplastic spring (the force with which the spring resists the elongation according to Hooke's law) doesn't longer increase when the spring is stretched beyond the critical value. The maximal stress the elastoplastic spring can achieve is called spring's elastic limit. Therefore, when the distance between two nodes of an elastoplastic  lattice spring model is locked and an is governed by a gradually increasing function, there will be a maximal stress that the combined system can develop in response to the stretching. This maximal stress is called {\it response force} $F$ and it is a measure of the strength of the material. Minimization of the cost of the material is usually carried out under assumption that the strength is not compromised, e.g. under condition
\begin{equation}\label{Fconstraint}
      F\ge const.
\end{equation}
A general procedure for computation of the response force of elastoplastic lattice-spring models is addressed in Gudoshnikov et al \cite{gudo} (which converts the problem to a so-called sweeping process linking optimization of the response force to the recent methods of e.g. Kamenskii et al \cite{obu}). A possible way to compute the cost of the material is by summing up the elastic limits $c_j$ of all springs, i.e. by taking
\begin{equation}\label{C4springs}
     C=\sum_{j=1}^m c_j,
\end{equation}
where $m\in\mathbb{N}$ is the total number of springs in the system. Therefore, the problem of minimization of the cost of the material is the search for a topology of a lattice spring model that minimizes quantity (\ref{C4springs}) under constraint (\ref{Fconstraint}). As a matter of fact, such an optimal topology consists of just 1 spring, i.e. the aforementioned optimization problem is trivial. That is why in this paper we propose to view springs of the lattice spring model current conducting and replace (\ref{Fconstraint}) by
\begin{equation}\label{FRconstraint}
      const\cdot F+const\cdot R\ge const,
\end{equation}
where $R$ is the electrical resistance of the lattice spring model viewed as an electric circuit. Design of multi-functional materials of this type is addressed in e.g. \cite{Torquato2} 
Because one possible realization of a current conducting spring can be achieved by fiber-reinforced polymer matrix composites, then following \cite{ch5mater01}, we consider resistance of each particular spring to be equal 
\begin{equation}
   R_i=1/c_i,
\label{ch5eq_Ri}  
\end{equation}
thus the total resistance of the system is a function of $c_1,...,c_m$ too.
The search for topology of the lattice spring model that minimizes the cost (\ref{C4springs}) under constraint (\ref{FRconstraint}) is the problem we address. 

\vskip0.2cm

\noindent Because numeric simulations of Malhotra et al \cite{arXiv} showed that optimal topologies of a network of 5 springs always attained at the values of parameters $c_1,...,c_5$ where at least one of $c_i$ is zero, the present paper  specifically focuses on lattice spring models of 4 springs. For instructional purposes, we will consider particular constants in (\ref{Fconstraint}) and (\ref{FRconstraint}) (to demonstrate computations explicitly), but arbitrary constants and general results can be obtained along the same lines.    

\vskip0.2cm

\noindent The paper is organized as follows. In the next section we introduce the details of fabrication cost optimization problem, which appear to yield 10 cases in total to be considered. The minimal fabrication cost for one  of these 10 cases is explicitly computed in Section~\ref{secCase9} and Section~\ref{secOthers} subsequently shows that the remaining cases cannot beat the cost obtained in Section~\ref{secCase9}.  
Concise summaries of computations for each of the cases is arrange in tables and the details of computations are sent to the Appendix. Wolfram Mathematica notebook used for graphic verification of the constraints is attached as Supplementary Material.

\section{Formulation of the optimization problem}

\noindent The list of all arrangements of 4 springs on 1-dimensional nodes is given in Table~\ref{list}. We will stick to the following particular version of (\ref{FRconstraint})
$$
  F_R(c)\ge 0.5,\quad F_R(c)=0.2\cdot F(c)+0.1\cdot R(c),\quad c=(c_1,c_2,c_3,c_4).
$$

\begin{table}[h]
\begin{center}   
\begin{tabular}{|c|c|}
  \hline
  \makecell{\includegraphics[width=0.4\textwidth]{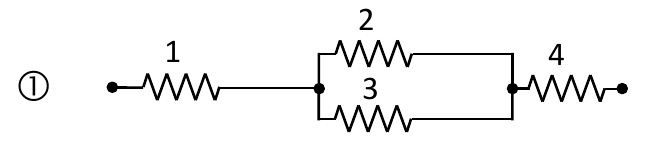}} & \makecell{\includegraphics[width=0.4\textwidth]{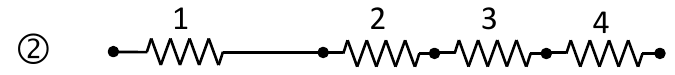}} \\
  \hline
  \makecell{\includegraphics[width=0.4\textwidth]{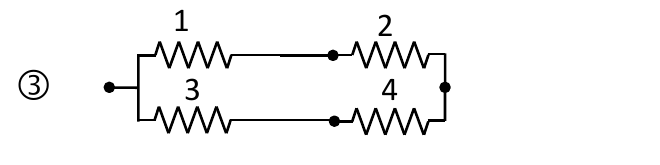}} & \makecell{\includegraphics[width=0.4\textwidth]{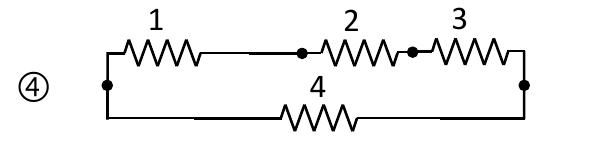}} \\
  \hline
  \makecell{\includegraphics[width=0.4\textwidth]{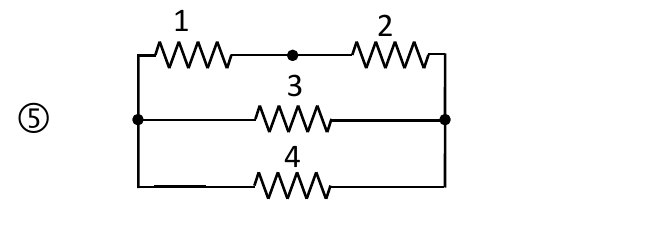}} & \makecell{\includegraphics[width=0.4\textwidth]{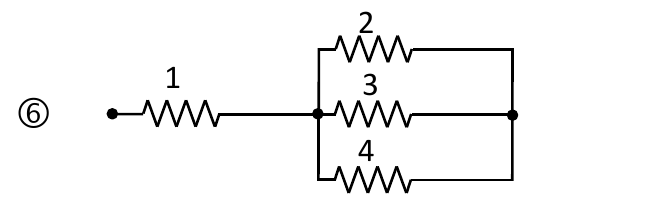}} \\
  \hline
  \makecell{\includegraphics[width=0.4\textwidth]{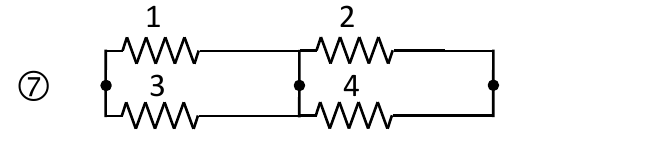}} & \makecell{\includegraphics[width=0.4\textwidth]{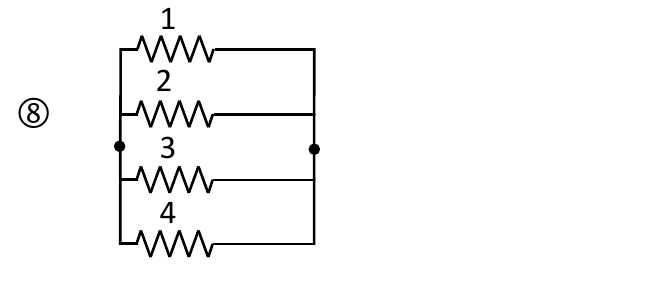}} \\
  \hline
  \makecell{\includegraphics[width=0.4\textwidth]{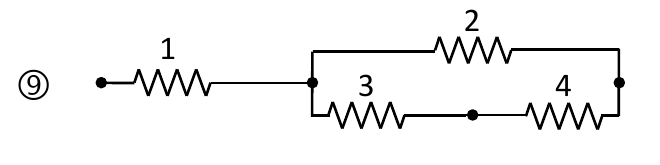}} & \makecell{\includegraphics[width=0.4\textwidth]{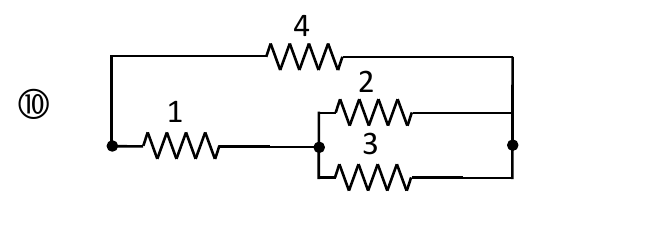}} \\
  \hline
\end{tabular}
\end{center}
\caption{List of all arrangements of 4 springs on 1-dimensional nodes. The arrangement with number $N$ in a circle is referred to as Case~$N$ in the sequel.}\label{list}
\end{table}

\noindent While minimizing the cost across the topologies from Table~\ref{list}, we don't want the response force to be compromised. For this reason we additionally consider the following particular version of (\ref{Fconstraint})
$
  F(c)\ge 0.75.
$ To conclude, our final optimization problem reads as
\begin{equation}\label{optp}
    F(c)\ge 0.75,\quad F_R(c)\ge 0.5,\quad C(c)\to \min.
\end{equation}


\section{Reduction to 2 variables and analytic computation of the minimal cost $C_{min}$ for Case 9.} \label{secCase9}

\subsection{Case 9}

\begin{center}   
\includegraphics[width=0.4\textwidth]{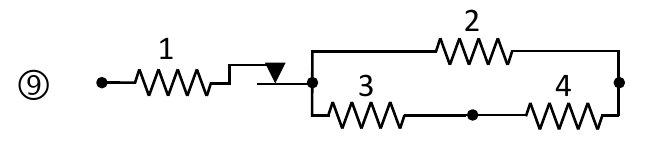}\qquad\quad\includegraphics[width=0.4\textwidth]{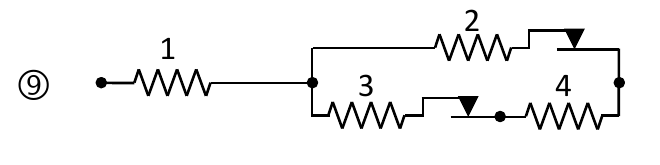}
\end{center}

\noindent{\bf Resistance.} Following standard rules of the computation of resistance of parallel and series connection of resistances \cite{Nilsson}, we compute
\begin{equation}\label{resistance9}
    R=R_1 + \frac{1}{\frac{1}{R_4}+\frac{1}{R_2+R_3}}
=\frac{1}{c_1} + \frac{1}{c_4 + \frac{1}{\frac{1}{c_2} + \frac{1}{c_3}}}.
\end{equation}

\noindent {\bf Response force.} To understand the total response force of the entire system, we can view the three springs $c_2$, $c_3$, $c_4$ as a single equivalent spring, whose elastic bound $c_{234}$ is to be determined. Then 
\begin{equation}\label{c1<}
     c_1<c_{234}
\end{equation}
would imply that, under gradual displacement-controlled loading, spring $c_1$ reaches elastic limit before $c_{234}$ does so, and the response force of the entire system cannot go beyond $c_1$. Likewise, when
\begin{equation}\label{c1>}
     c_1>c_{234},
\end{equation}
the response force of the system is $c_{234}$. 

\vskip0.2cm

\noindent Determining $c_{234}$ means understanding the maximal response force of the network of springs $c_2,$ $c_3,$ $c_4.$ Here we can again replace the series connection of $c_2$ and $c_3$ by an equivalent spring whose elastic bound $c_{23}$ still needs to be determined. If
\begin{equation}\label{c_2<c_3}
      c_2<c_3
\end{equation}
then, under gradual stretching of the series connection of springs $c_2$ and $c_3$, the response force of this connection stops increasing after reaching the value $c_{23}=c_2$ (beginning this instance, spring $c_2$ will just stretch plastically and keep the same stress $c_2$). As for spring $c_4$, it is maximal stress is $c_4$. We recall that the total force required to move the object a given displacement will be equal to the sum of the forces required to move each individual spring by the same amount. Therefore, the maximal stress of the network of springs $c_2$, $c_3$, $c_4$ will be
\begin{equation}\label{c234a}
  c_{234}=c_2+c_4
\end{equation}
if condition (\ref{c_2<c_3}) holds. Analogously,
if 
$$
      c_2>c_3
$$
then
\begin{equation}\label{c234b}
    c_{234}=c_3+c_4.
\end{equation}
Substituting (\ref{c234a}) and (\ref{c234b}) to (\ref{c1<}) and (\ref{c1>}), we get 4 cases
$$
\begin{array}{ll}
   c_2<c_3,\ c_1<c_2+c_4: & F=c_1,\\
   c_2<c_3,\ c_1>c_2+c_4: & F=c_2+c_4,\\
   c_2>c_3,\ c_1<c_3+c_4: & F=c_1,\\
   c_2>c_3,\ c_1>c_3+c_4: & F=c_3+c_4.
\end{array}
$$
Out of these 4 cases those with $c_2<c_3$ and those with $c_2>c_3$ are equivalent in terms of further analysis, so we won't restrict the generality by focusing on just $c_2<c_3,$ whose derivation can be  formulated through Table~\ref{case9table}.

\vskip0.2cm

\begin{table}[h]
\begin{tabular}{l|l}
Response force $F$ & Multi-functional performance $F_R$ and cost $C$\\ \hline
$\begin{array}{c}
c_2<c_3\\
c_{234}=c_2+c_4
\end{array}$              &   $\begin{array}{l}F_R(c)=0.2F(c)+0.1R(c),\quad {\rm where}\ c=(c_1,c_2,c_3,c_4)\\
C(c)=c_1+c_2+c_3+c_4\end{array}$                           \\ \hline
$\begin{array}{c}
\mbox{\bf Case 9.1}\\
c_1<c_2+c_4\\
F=c_1 \\ c_1\mbox{ is plastic}
\end{array}$                 & $
\begin{array}{l}\tilde F_R(c_1,c_2):=0.2c_1 + 0.1 \left(\frac{1}{c_1}+\frac{1}{c_1-\frac{c_2}{2}}\right)\ge F_R(c)\\
\tilde C(c_1,c_2):=2c_1+c_2\le C(c)\end{array}$                            \\ \hline
$\begin{array}{c}
\mbox{\bf Case 9.2}\\
c_1>c_2+c_4\\
F=c_2+c_4\\
c_2\mbox{ is plastic}
\end{array}$              & $
\begin{array}{l}\tilde F_R(c_2,c_4):=0.2(c_2 + c_4) + 0.1\cdot\left(\frac{1}{c_2+c_4}+\frac{1}{\frac{c_2}{2} + c_4}\right)\ge F_R(c)\\
\tilde C(c_2,c_4):=3c_2+2c_4\le C(c)\end{array} $                     
\end{tabular}
\caption{
 The two different plastic behaviors possible in the case~9 (left column) and the corresponding two-variable estimate for the $F_R$ and $C$ (right column). The derivations of the estimates are given in Appendix. The constant $c_{234}$ refers to the strength of the the part of the network formed by just springs $c_1$, $c_3$, $c_4$ when both springs $c_2$ and $c_4$ reached their maximal strength.}
\label{case9table}
\end{table}

\noindent Optimization problem (\ref{optp}) corresponding to Case 9.1 will be taken as follows
\begin{eqnarray}\label{opt9}
\hskip-2cm &&\mbox{\bf Case 9.1:} \ c_2\le c_3,\ c_1\le c_2+c_4, \ c_1\ge 0.75,\  F_R(c)\ge 0.5,\  C(c)\to {\rm min},
\end{eqnarray}
considered in the positive quadrant $c_i\ge 0,$ $i=1,2,3,4.$
Along with full problem (\ref{opt9reduced}) we consider the following reduced problem with notations from Table~\ref{case9table}:
\begin{eqnarray}\label{opt9reduced}
\hskip-2cm &&\mbox{\bf Case 9.1 reduced:}
 \ c_1\ge 0.75,\ c_2\ge 0,\ \tilde F_R(c_1,c_2)\ge 0.5,\ \tilde C(c_1,c_2)\to {\rm min}
\end{eqnarray}

\noindent The inequalities in (\ref{opt9}) and (\ref{opt9reduced}) will be referred to as full domain and reduced domain accordingly.

\begin{proposition} \label{reducedproposition}
    If $(c_{1*},c_{2*})$ solves reduced problem (\ref{opt9reduced}) then
\begin{equation}\label{c*}
      c_*=(c_{1*}, c_{2*}, c_{2*}, c_{1*}- c_{2*}) 
\end{equation}
solves full problem (\ref{opt9}).
\end{proposition}

\noindent {\bf Proof.} Since $( c_{1*}, c_{2*})$ belongs to the reduced domain, $c_*$ belongs to the full domain.
Since the minimum of $\tilde C$ is attained at $(c_{1*}, c_{2*})$, 
$$
  \tilde C(c_{1*}, c_{2*})\le \tilde C(c_1,c_2),
$$
for any $(c_1,c_2)$ from the reduced domain. 
Secondly,
$$
  \tilde C(c_1,c_2)\le C(c_1,c_2,c_3,c_4)
$$
for any $(c_1,c_2,c_3,c_4)$ from the full domain. Finally, 
$$
   \tilde C(c_{1*},c_{2*})=C(c_*)
$$
by direct inspection. To summarize, we proved that
$$
  C(c_*)\le C(c_1,c_2,c_3,c_4) 
$$
for all $(c_1,c_2,c_3,c_4)$ from the full domain. \qed

\vskip0.5cm

\begin{figure}[h]
\begin{center}
    \includegraphics[scale=0.8]{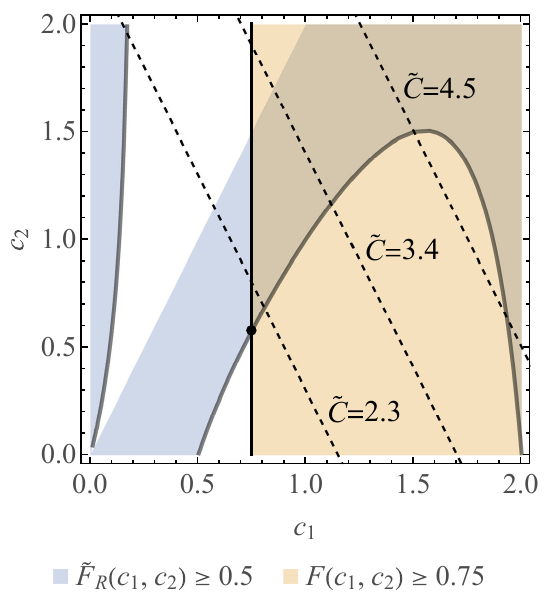}
    \includegraphics[scale=0.786]{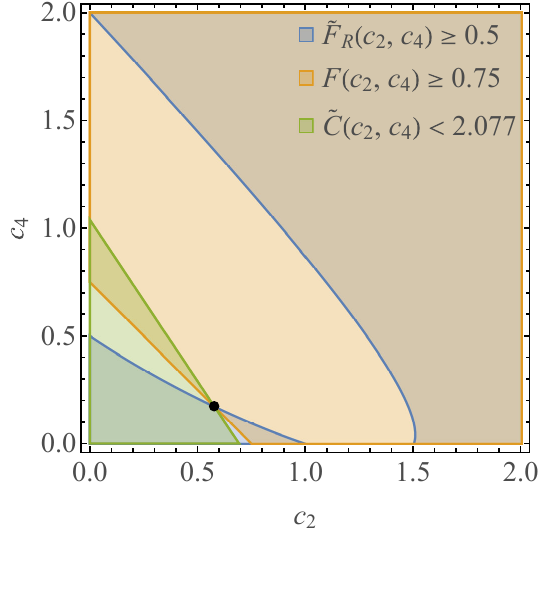}
\caption{Left: Visualization of the two-variable estimate for Case 9.1, Right: Two-variable estimate for Case 9.2.}
\label{samples}
\end{center}

\end{figure}

\noindent It turns of that the reduced optimization problem (\ref{opt9reduced}) can be solved in an elementary way. Indeed, the graph of $\tilde C(c_1,c_2)=const$ is a line, whose samples for different values of the constant are shown in Fig.~\ref{samples}\,(left). From Fig.~\ref{samples} we conclude that in order to solve problem (\ref{opt9}) (whose solution is denoted by a bold dot in Fig.~\ref{samples}\,(left)), we have to find the intersection (marked with bold dot) of the leftmost boundary of the region $F(c_1)\ge 0.75)$ with the bottom boundary of $\tilde F_R(c_1,c_2)\ge 0.5$. As the graph illustrates, this point of intersection is the intersection of $c_1=0.75$ with $\tilde F(0.75,c_2)=0.5$, i.e. solution of 
$$
0.2\cdot0.75+0.1\left(\frac{1}{0.75}+\frac{1}{0.75-\frac{c_2}{2}}\right)=0.5,
$$
which is
$
  c_2\approx 0.577.
$
And using Proposition~\ref{reducedproposition}, the solution of problem (\ref{opt9}) is
\begin{eqnarray}\label{opt9solution}
\hskip-2cm &&\mbox{\bf Case 9.1  solution:}
 \quad c_{1*}=0.75,\quad c_{2*}\approx 0.577,\quad c_{3*}\approx 0.577,\quad c_{4*}\approx 0.173
\end{eqnarray}
In particular, the minimum of $\tilde C(c_1,c_2)$ in problem (\ref{opt9reduced}) computes as
\begin{equation}\label{mincost}
      \tilde C(0.75,0.577)\approx 2.077.
\end{equation}    
What we show in the remainder of the paper is that none of the other configurations from the list of Table~\ref{list} are capable to produce the cost smaller than (\ref{mincost}).

\vskip0.2cm

\noindent Because the optimal value (\ref{opt9solution}) belongs to the boundary between the cases 9.1 and 9.2, it is expected that (\ref{opt9solution}) will also be a solution to
\begin{eqnarray}\label{opt92}
\hskip-2cm &&\mbox{\bf Case 9.2:} \ c_2\le c_3,\ c_1\ge c_2+c_4, \ c_2+c_4\ge 0.75,\ F_R(c)\ge 0.5,\  C(c)\to {\rm min},
\end{eqnarray}
considered in the positive quadrant $c_i\ge 0,$ $i=1,2,3,4.$ Indeed,
if we consider the corresponding reduced problem 
\begin{eqnarray}\label{opt92reduced}
\hskip-2cm &&\mbox{\bf Case 9.2 reduced:}
 \ c_2+c_4\ge 0.75,\ c_2\ge 0,\ \tilde F_R(c_1,c_2)\ge 0.5,\ \tilde C(c_1,c_2)\to {\rm min},
\end{eqnarray}
Fig.~\ref{samples}\,(right)  shows the cost  (\ref{mincost}) is possible for a single pair of parameters  $(c_2,c_4)$ coinciding with the values in (\ref{opt9solution}). Fig.~\ref{samples}\,(right) also impliesj that none of $(c_2,c_4)$ with $\tilde C(c_2,c_4)<2.077$ intersect the reduced domain of (\ref{opt92reduced}), so the cost smaller than (\ref{mincost}) is impossible for the full problem (\ref{opt92}) either as follows rigorously from the next proposition.

\vskip0.2cm

\begin{proposition}\label{prop2}
If $(c_1,c_{2},c_{3},c_{4})$ belongs to the full domain, then $(c_{1},c_{2})$ belongs to the reduced domain and 
\begin{equation}\label{less}    
   \tilde C(c_{1},c_{2})\le C(c_{1},c_{2},c_{3},c_{4}).
\end{equation}
In particular, if the reduced domain doesn't have $(c_1,c_2)$ satisfying 
$$
  \tilde C(c_1,c_2)\le C_*,
$$
then minimum of the full problem is greater than $C_*$.
\end{proposition}

\noindent {\bf Proof.} Indeed, in the settings of proposition we have
$$
  C_*<\tilde C(c_1,c_2)
$$
for all $(c_1,c_2)$ from the reduced domain. Substituting (\ref{less}), we conclude
$$
  C_*<C(c_{1},c_{2},c_{3},c_{4})
$$
for all $(c_1,c_2,c_3,c_4)$ from the full domain.\qed

\noindent Now we proceed with the rest of the cases of Table~\ref{list} which will all have a cost greater than (\ref{mincost}) by \ref{prop2}.

\section{Reduction to 1 or 2 variables and non-reachability of the cost of $C_{min}$ for all other cases.}\label{secOthers}

\subsection{Case 1.} \label{secCase1} The resistance in this case computes as


\begin{equation}\label{opt1}
R=R_1+\frac{1}{\frac{1}{R_2}+\frac{1}{R_3}}+R_4=\frac{1}{c_1}+\frac{1}{c_2+c_3}+\frac{1}{c_4}.
\end{equation}

\begin{table}[h]
\begin{tabular}{ll}
\multicolumn{1}{l|}{Response force $F$} & Multi-functional performance $F_R$ and cost $C$\\ \hline
$\begin{array}{c}
c_1<c_4\\
c_{23}=c_2+c_3
\end{array}$              &   $\leftarrow\begin{array}{l}
\mbox{when both }c_2\mbox{ and }c_3\mbox{ reach elastic bounds, the stress} \\
\mbox{of the component formed by springs }c_2\mbox{ and }c_3\mbox{ is }c_2+c_3\end{array}$                           \\ \hline
\multicolumn{1}{l|}{$\begin{array}{l}
\mbox{\bf Case 1.1}\\
c_1<c_2+c_3\\
F=c_1 \\ c_1\mbox{ is plastic}
\end{array}$                 }& $  \def\arraystretch{1.5}
\begin{array}{l}\tilde F_R(c_1):=0.2c_1 + 0.1 \cdot\frac{3}{c_1}\ge F_R(c)\\
\tilde C(c_1):=3c_1\le C(c)\end{array}$                            \\ \hline
\multicolumn{1}{l|}{$\begin{array}{l}
\mbox{\bf Case 1.2}\\
c_1>c_2+c_3\\
F=c_2+c_3\\
c_2,c_3\mbox{ are plastic}
\end{array}$              }& $\def\arraystretch{1.5}
\begin{array}{l}\tilde F_R(c_2+c_3):=0.2(c_2 + c_3) + 0.1\cdot \frac{3}{c_2+c_3}\ge F_R(c)\\
\tilde C(c_2+c_3):=3(c_2+c_3)\le C(c)\end{array} $                     
\end{tabular}
\caption{Summary of a one-variable estimate for $F_R$ and $C$ in the Case 1, see Appendix (Section~\ref{SectionApp}) for the details of computations.}
\label{case1table}
\end{table}

\begin{figure}[h]
    \centering
    \vskip0.2cm
    \includegraphics[width=0.8\textwidth]{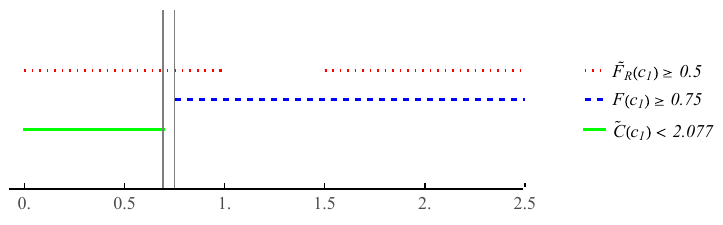}
    
    \caption{Visualization of the single variable estimate for Case 1.}
    \label{case1fig}
\end{figure}

\noindent Since formulas of Table~\ref{case1table} for Cases 1.1 and 1.2 are equivalent, it is sufficient to address just Case 1.1. Referring to $F$, $\tilde F_R$ and $\tilde C$ from Case 1.1 of Table~\ref{case1table}, the graph of Fig.~\ref{case1fig} shows the lack of the intersection of the 3 regions from the statement of Proposition~\ref{prop2}, which implies that the cost in the Case 1.1 is greater than 2.077 (the optimal value obtained in Case 9). 

\subsection{Case 2.} The resistance computes as

\begin{equation}\label{opt2}
R= R_1+R_2+R_3+R_4 = \frac{1}{c_1}+\frac{1}{c_2}+\frac{1}{c_3}+\frac{1}{c_4}.
\end{equation}

\begin{figure}[h]
    \centering
    \vskip0.2cm
    \includegraphics[width=0.8\textwidth]{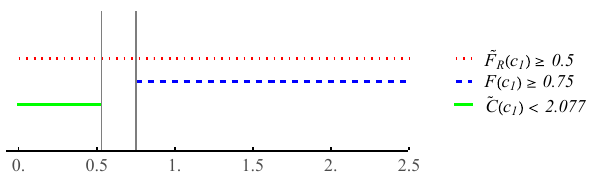}
    
    \caption{Vizualization of the single variable estimate for Case 2.}
    \label{case2fig}
\end{figure}

\noindent The values of $F$, $\tilde F_R$ and $\tilde C$ for the reduced optimization problem are summarized in Table~\ref{case2table} and visualized in Fig.~\ref{case2fig}, which, similarly to Case 1 (Section~\ref{secCase1}), confirms that the cost in Case 2 is higher than 2.077.

\vskip0.2cm

\begin{table}[h]
\begin{tabular}{ll}
\multicolumn{1}{l|}{Response force $F$} & Multi-functional performance $F_R$ and cost $C$ \\ \hline
$\begin{array}{c}
c_1 < c_2 ,  c_1 < c_3 , c_1 < c_4\\

\end{array}$ &
$\leftarrow\begin{array}{l}
{\rm When  } c_1  \mbox{ reaches elastic bound, the stress} \\
\mbox{of the component formed by springs } c_1,
\ c_2,\ c_3,\ c_4 \mbox{ is } c_1
\end{array}$ \\ \hline
\multicolumn{1}{l|}{$\begin{array}{l}
\mbox{\textbf{Case 2}} \\

F = c_1 \\
c_1  \mbox{ is plastic}
\end{array}$} &
$\def\arraystretch{1.5}
\begin{array}{l}
\tilde{F}_R(c_1) := 0.2(c_1) + 0.1 \cdot \frac{4}{c_1} \geq F_R(c) \\
\tilde{C}(c_1) := 4c_1  \leq C(c)
\end{array}$ \\ \hline
\end{tabular}
\caption{Summary of a one-variable estimate for $F_R$ and $C$ in the Case 2. See Section~\ref{SectionApp} for the details of computations.}
\label{case2table}
\end{table}

\subsection{Case 3.} The resistance is now

\begin{equation}\label{opt3}
R= \frac{1}{\frac{1}{R_1+R_2}+\frac{1}{R_3+R_4}}=\frac{1}{\frac{1}{\frac{1}{c_1}+\frac{1}{c_2}}+\frac{1}{\frac{1}{c_3}+\frac{1}{c_4}}}.
\end{equation}

\begin{figure}[h]
    \centering
    \vskip0.2cm
    \includegraphics[scale=1.3,trim={1.5cm 23cm 0 2.5cm}]{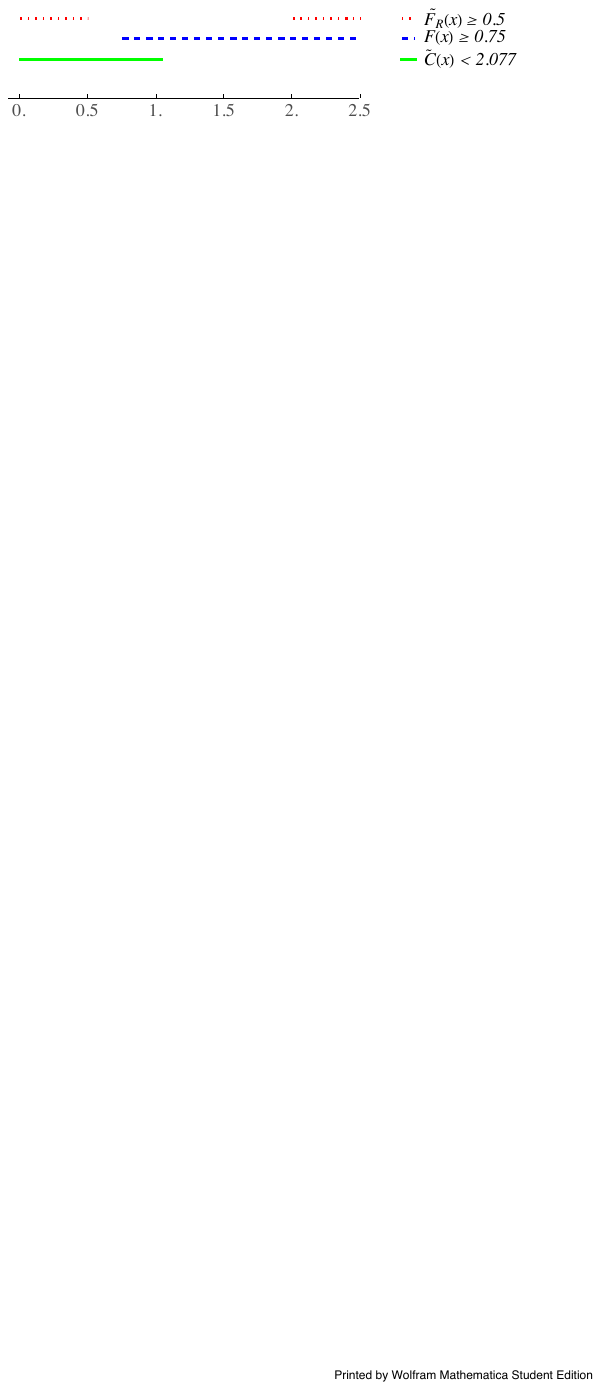}
    
    \caption{Vizualization of the single variable estimate for Case 3 with $x=c_1+c_3$.}
    \label{case3fig}
\end{figure}

\noindent The reader is referred to Appendix (Section~\ref{SectionApp}) for verification of the estimates of Table~\ref{case3table}. 
\noindent Based on Proposition~\ref{prop2} and Figures~\ref{case3fig} we conclude that (similarly to Case 4) the minimum of the full optimization problem corresponding to Case~3 is greater than the cost (\ref{mincost}) obtained for Case 9.

\vskip0.2cm

\begin{table}[h]
\begin{tabular}{ll}
\multicolumn{1}{l|}{Response force $F$} & Multi-functional performance $F_R$ and cost $C$ \\ \hline
$\begin{array}{c}
c_1 < c_2 ,  c_3 < c_4 \\

\end{array}$ &
$\leftarrow\begin{array}{l}
\mbox{When  } c_1 \mbox{ and } c_3 \mbox{ reaches elastic bound, the stress} \\
\mbox{of the component formed by springs } c_1,
\ c_2, \mbox{ is } c_1 \mbox{ and } c_3,\  c_4 \mbox{ is  } c_3
\end{array}$ \\ \hline
\multicolumn{1}{l|}{$\begin{array}{l}
\mbox{\textbf{Case 3}} \\

F = c_1 + c_3\\
c_1,\ c_3 \mbox{ are plastic}
\end{array}$} &
$\def\arraystretch{1.5}
\begin{array}{l}
\tilde{F}_R(c_1+c_3) := 0.2(c_1+c_3) + 0.1 \cdot \frac{2}{c_1+c_3} \geq F_R(c) \\
\tilde{C}(c_1+ c_3) := 2c_1+2c_3  \leq C(c)
\end{array}$ \\ \hline
\end{tabular}
\caption{Summary of a one-variable estimate for $F_R$ and $C$ in the Case 3.}
\label{case3table}
\end{table}

\FloatBarrier


\subsection{Case 4.} 
The resistance is 
\begin{equation}\label{opt4}
R=\frac{1}{R_1+R_2+R_3}+\frac{1}{R_4}=\frac{1}{\frac{1}{\frac{1}{c_1}+\frac{1}{c_2}+\frac{1}{c_3}}+c_4}
\end{equation}


\begin{table}[h]
\begin{tabular}{ll}
\multicolumn{1}{l|}{Response force $F$} & Multi-functional performance $F_R$ and cost $C$ \\ \hline
$\begin{array}{c}
c_1 < c_2 ,  c_1 < c_3 \\
c_{123} = c_1
\end{array}$ &
$\leftarrow\begin{array}{l}
\mbox{When  } c_1  \mbox{ reaches elastic bound, the stress} \\
\mbox{of the component formed by springs } c_1,
\ c_2,\ c_3 \mbox{ is } c_1
\end{array}$ \\ \hline
\multicolumn{1}{l|}{$\begin{array}{l}
\mbox{\textbf{Case 4}} \\

F = c_1+c_4 \\
c_1 \mbox{ and } c_4 \mbox{ are plastic}
\end{array}$} &
$\def\arraystretch{1.5}
\begin{array}{l}
\tilde{F}_R(c_1,\ c_4) := 0.2(c_1+c_4) + 0.1 \cdot \frac{1}{\frac{c_1}{3}+c_4} \geq F_R(c) \\
\tilde{C}(c_1,\ c_4) := 3c_1 + c_4 \leq C(c)
\end{array}$ \\ \hline
\end{tabular}
\caption{Summary of a two-variable estimate for $F_R$ and $C$ in the Case 4, see  Appendix (Section~\ref{SectionApp}) for details. 
}
\label{case4table}
\end{table}

\begin{figure}[h]
\begin{center}
    \includegraphics[scale=0.773]{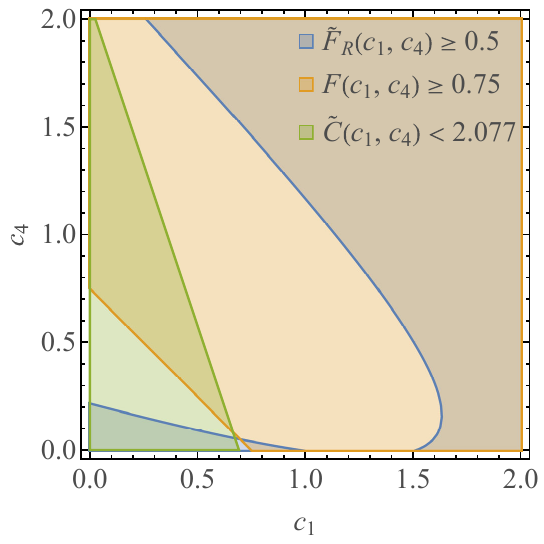}
    \includegraphics[scale=0.8]{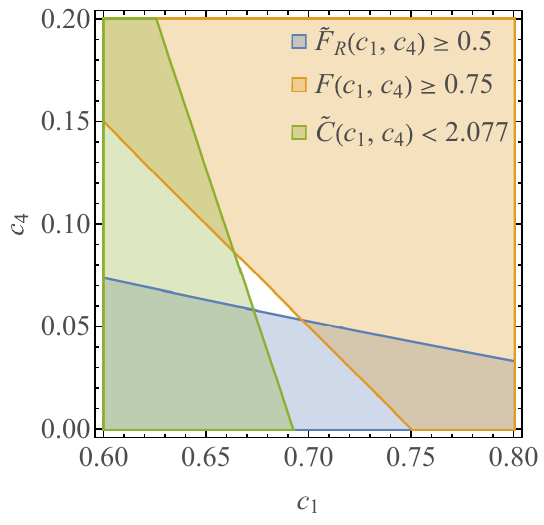}
\caption{Left: Visualization of the two-variable estimates for Case 4, Right: Left Figure zoomed in.}
\label{figcase4}
\end{center}
\end{figure}

\noindent Based on Proposition~\ref{prop2} and Figures~\ref{figcase4} (that shows lack of intersection of the 3 regions) we use Proposition~\ref{prop2} to conclude that the minimum of the full optimization problem corresponding to Case~4 is greater than the cost (\ref{mincost}) obtained for 
Case 9.
\FloatBarrier





\subsection{Case 5.} The resistance is 
\begin{equation}\label{opt5}
R=\frac{1}{\frac{1}{R_1+R_2}+\frac{1}{R_3}+\frac{1}{R_4}} = \frac{1}{\frac{1}{\frac{1}{c_1}+\frac{1}{c_2}}+c_3+c_4}
\end{equation}

\begin{table}[h]
\begin{tabular}{ll}
\multicolumn{1}{l|}{Response force $F$} & Multi-functional performance $F_R$ and cost $C$ \\ \hline
$\begin{array}{c}
c_1 < c_2  \\

\end{array}$ &
$\leftarrow\begin{array}{l}
\mbox{When  } c_1  \mbox{ reaches elastic bound, the stress} \\
\mbox{of the component formed by springs } c_1,
\ c_2, \mbox{ is } c_1 
\end{array}$ \\ \hline
\multicolumn{1}{l|}{$\begin{array}{l}
\mbox{\textbf{Case 5}} \\

F = c_1 + c_3 + c_4\\
c_1,\ c_3,\ c_4 \mbox{ are plastic}
\end{array}$} &
$\def\arraystretch{1.5}
\begin{array}{l}
\tilde{F}_R(c_1,c_3+c_4) := 0.2(c_1+c_3+c_4) + 0.1 \cdot \frac{1}{\frac{1}{\frac{2}{c_1}}+c_3+c_4} \geq F_R(c) \\
\tilde{C}(c_1, c_3+c_4) := 2c_1+c_3+c_4  \leq C(c)
\end{array}$ \\ \hline
\end{tabular}
\caption{Summary of two-variable estimate for $F_R$ and $C$ in the Case 5,  see  Appendix (Section~\ref{SectionApp}) for details.}
\label{case5table}
\end{table}

\begin{figure}[h]
    \centering
    \vskip0.2cm
    \includegraphics[scale=0.9, trim={0cm 16cm 5cm 3cm}]{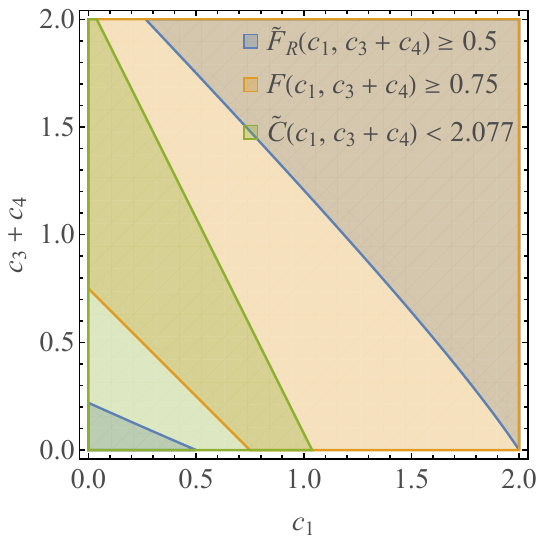}
    \caption{Visualization of the two-variable estimates for Case 5.}
    \label{figcase5}
\end{figure}

\noindent For formulas of Table~\ref{case5table} for visualization $\tilde{F}_R$ and $\tilde{C}$ can be expressed in 2 variables, taking $x=c_1$ and $y=c_3+c_4$, see Fig.~\ref{figcase5}, which shows lack of intersection of the 3 regions.

\subsection{Case 6.} 


\begin{eqnarray}\label{opt6}
\hskip-2.32cm &&\mbox{Resistance:} \ R=R_1+\frac{1}{\frac{1}{R_2}+\frac{1}{R_3}+\frac{1}{R_4}}=\frac{1}{c_1}+\frac{1}{c_2+c_3+c_4}
\end{eqnarray}


\vskip0.2cm

\begin{table}[h]
\begin{tabular}{ll}
\multicolumn{1}{l|}{Response force $F$} & Multi-functional performance $F_R$ and cost $C$\\ \hline
$\begin{array}{c}
c_{234}=c_2+c_3+c_4
\end{array}$              &   $\leftarrow\begin{array}{l}
\mbox{when }c_2,\ c_3, \mbox{ and } c_4\mbox{ reach elastic bounds, the stress} \\
\mbox{of the component formed by springs }c_2,\ c_3 \mbox{ and }c_4\mbox{ is }c_2+c_3+c_4\end{array}$                           \\ \hline
\multicolumn{1}{l|}{$\begin{array}{l}
\mbox{\bf Case 6.1}\\
c_1<c_2+c_3+c_4\\
F=c_1 \\ c_1\mbox{ is plastic}
\end{array}$                 }& $  \def\arraystretch{1.5}
\begin{array}{l}\tilde F_R(c_1):=0.2c_1 + 0.1 \cdot\frac{2}{c_1}\ge F_R(c)\\
\tilde C(c_1):=2c_1\le C(c)\end{array}$                            \\ \hline
\multicolumn{1}{l|}{$\begin{array}{l}
\mbox{\bf Case 6.2}\\
c_1>c_2+c_3+c_4\\
F=c_2+c_3+c_4\\
c_2,c_3,c_4\mbox{ are plastic}
\end{array}$              }& $\def\arraystretch{1.5}
\begin{array}{l}\tilde F_R(c_2+c_3+c_4):=0.2(c_2 + c_3 + c_4) + \frac{0.2}{c_2+c_3+c_4}\ge F_R(c)\\
\tilde C(c_2+c_3+c_4):=2(c_2+c_3+c_4)\le C(c)\end{array} $                     
\end{tabular}
\caption{Summary of a one-variable estimate for $F_R$ and $C$ in the Case 6.}
\label{case6table}
\end{table}

\begin{figure}[h]
    \centering
    \vskip0.2cm
    \includegraphics[width=0.8\textwidth]{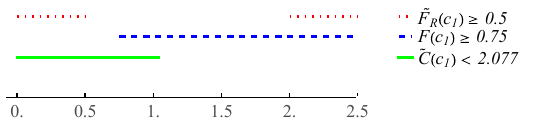}
    
    \caption{Visualization of the single variable estimate for Case 6.}
    \label{case6fig}
\end{figure}

\noindent Since formulas of Table~\ref{case6table} for Cases 6.1 and 6.2 are equivalent, it is sufficient to address just Case 6.1. And based on Figure~\ref{case6fig} we conclude that minimal cost in Case 6 cannot get lower than the one obtained for Case 9 either.

\subsection{Case 7.} 

\noindent The resistance is 
\begin{equation}\label{opt7}
R= \frac{1}{\frac{1}{R_1}+\frac{1}{R_3}}+\frac{1}{\frac{1}{R_2}+\frac{1}{R_4}}=\frac{1}{c_1+c_3}+\frac{1}{c_2+c_4}
\end{equation}


\vskip0.2cm

\begin{table}[h]
\begin{tabular}{ll}
\multicolumn{1}{l|}{Response force $F$} & Multi-functional performance $F_R$ and cost $C$\\ \hline
$\begin{array}{c}
c_{13}=c_1+c_3\\
c_{24}=c_2+c_4
\end{array}$              &   $\leftarrow\begin{array}{l}
\mbox{when }c_1\mbox{ and } c_3 \mbox{ or } c_2\mbox{ and } c_4 \mbox{ reach elastic bounds, the stress} \\
\mbox{of the component formed by springs }c_1 \mbox{ and }c_3\mbox{ is }c_1+c_3 \mbox{ and } \mbox{the }\\
\mbox{component formed by springs }c_2 \mbox{ and }c_4\mbox{ is }c_2+c_4\end{array}$                           \\ \hline
\multicolumn{1}{l|}{$\begin{array}{l}
\mbox{\bf Case 7.1}\\
c_1+c_3<c_2+c_4\\
F=c_1+c_3 \\ c_1,c_3\mbox{ are plastic}
\end{array}$                 }& $  \def\arraystretch{1.5}
\begin{array}{l}\tilde F_R(c_1+c_3):=0.2(c_1+c_3) + 0.1 \cdot\frac{2}{c_1+c_3}\ge F_R(c)\\
\tilde C(c_1+c_3):=2(c_1+c_3)\le C(c)\end{array}$                            \\ \hline
\multicolumn{1}{l|}{$\begin{array}{l}
\mbox{\bf Case 7.2}\\
c_1+c_3>c_2+c_4\\
F=c_2+c_4\\
c_2,c_4\mbox{ are plastic}
\end{array}$              }& $\def\arraystretch{1.5}
\begin{array}{l}\tilde F_R(c_2+c_4):=0.2(c_2 + c_4) + \frac{0.2}{c_2+c_4}\ge F_R(c)\\
\tilde C(c_2+c_4):=2(c_2+c_4)\le C(c)\end{array} $                     
\end{tabular}
\caption{Summary of a one-variable estimate for $F_R$ and $C$ in the Case 7, see Appendix (Section~\ref{SectionApp}) for details. }
\label{case7table}
\end{table}

\begin{figure}[h]
    \centering
    \vskip0.2cm
    \includegraphics[scale=1.3,trim={1.5cm 22.5cm 0 2.5cm}]{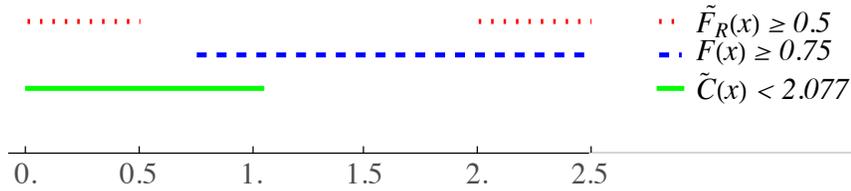}
    
    \caption{Visualization of the single variable estimate for Case 7 where $x=c_1+c_3$.}
    \label{case7fig}
\end{figure}

\noindent Since formulas of Table~\ref{case7table} for Cases 7.1 and 7.2 are equivalent, it is sufficient to address just Case 7.1, which again shows the cost greater than 2.077, see Fig.~\ref{case7fig}.

\subsection{Case 8.} 

\noindent The resistance:

\begin{equation}\label{opt8}
R=\frac{1}{\frac{1}{R_1}+\frac{1}{R_2}+\frac{1}{R_3}+\frac{1}{R_4}}=\frac{1}{c_1+c_2+c_3+c_4}
\end{equation}

\begin{figure}[h]
    \centering
    \vskip0.2cm
    \includegraphics[scale=1.4,trim={1.5cm 23cm 0 2.5cm}]{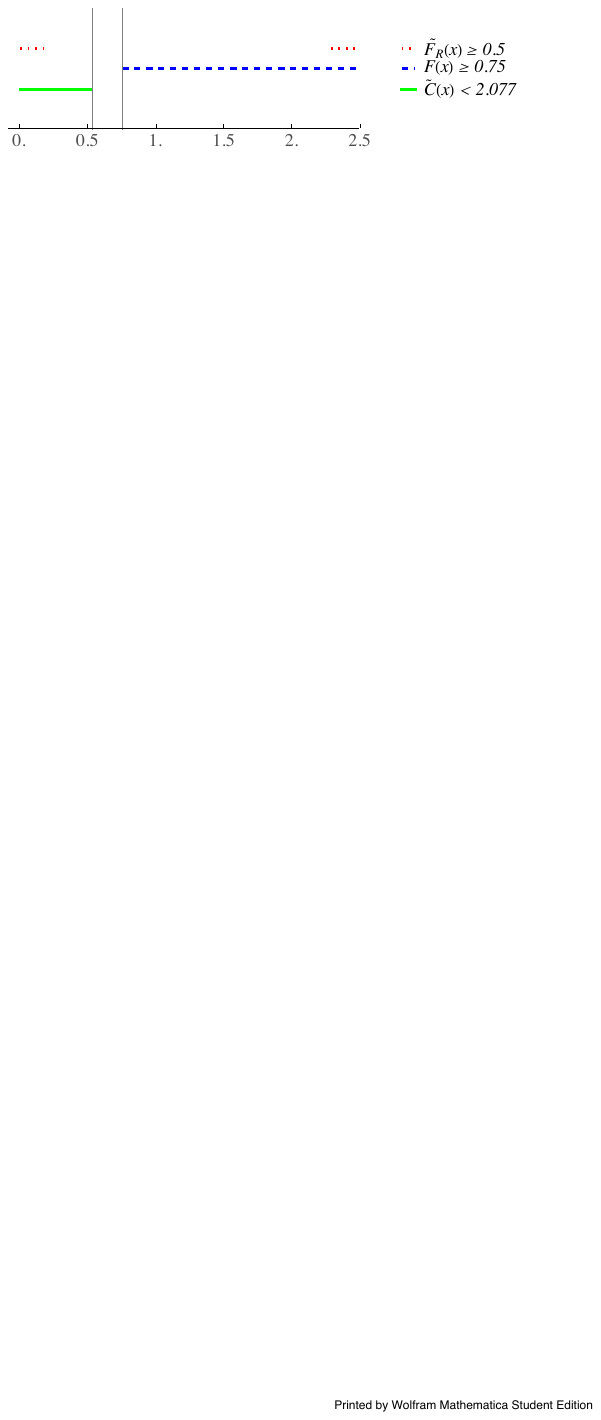}
    
    \caption{Vizualization of the single variable estimate for case 8 where $x=c_1+c_2+c_3+c_4$.}
    \label{case8fig}
\end{figure}

\noindent Based on Proposition~\ref{prop2} and Figure~\ref{case8fig} (that show lack of intersection) we conclude that the minimum of the full optimization problem corresponding to Case~8 is greater than the cost (\ref{mincost}) obtained for 
Case 9.

\vskip0.2cm

\begin{table}[h]
\begin{tabular}{ll}
\multicolumn{1}{l|}{Response force $F$} & Multi-functional performance $F_R$ and cost $C$ \\ \hline

\multicolumn{1}{l|}{$\begin{array}{l}
\mbox{\textbf{Case 8}} \\

F = c_1 +c_2+ c_3 + c_4\\
c_1,\ c_2,\ c_3,\ c_4 \mbox{ are plastic}
\end{array}$} &
$\def\arraystretch{1.5}
\begin{array}{l}
\tilde{F}_R(c_1+c_2+c_3+c_4) := 0.2(c_1+c_2+c_3+c_4) + 0.1 \cdot \\ \quad \quad \quad \quad \quad \quad \quad \quad \quad \quad   \quad \frac{1}{(c_1+c_2+c_3+c_4)}\geq F_R(c) \\
\tilde{C}(c_1+c_2 +c_3+c_4) := c_1+c_2+c_3+c_4  \leq C(c)
\end{array}$ \\ \hline
\end{tabular}
\caption{Summary of a one-variable estimate for $F_R$ and $C$ in the Case 8, see Appendix (Section~\ref{SectionApp}) for details. \\
}
\label{case8table}
\end{table}

\noindent For formulas of Table~\ref{case8table} for visualization $\tilde{F}_R$ and $\tilde{C}$ can be expressed in 1 variables, taking $x=c_1+c_2+c_3+c_4$ .

\subsection{Case 10.} 

\noindent {The resistance in this case computes as} 

\begin{equation}\label{opt10}
R= \frac{1}{\frac{1}{R_1+\frac{1}{\frac{1}{R_2}+\frac{1}{R_3}}}+\frac{1}{R_4}}=\frac{1}{c_4+\frac{1}{\frac{1}{c_1}+\frac{1}{c_2+c_3}}}
\end{equation}

\begin{table}[h]
\begin{tabular}{ll}
\multicolumn{1}{l|}{Response force $F$} & Multi-functional performance $F_R$ and cost $C$\\ \hline
$\begin{array}{c}
c_{23}=c_2+c_3
\end{array}$              &   $\leftarrow\begin{array}{l}
\mbox{when both }c_2\mbox{ and }c_3\mbox{ reach elastic bounds, the stress} \\
\mbox{of the component formed by springs }c_2\mbox{ and }c_3\mbox{ is }c_2+c_3\end{array}$                    \\ \hline
\multicolumn{1}{l|}{$\begin{array}{l}
\mbox{\bf Case 10.1}\\
c_1<c_2+c_3\\
F=c_1+c_4 \\ c_1,c_4\mbox{ are plastic}
\end{array}$                 }& $  \def\arraystretch{1.5}
\begin{array}{l}\tilde F_R(c_1,c_4):=0.2(c_1+c_4) + 0.1 \cdot\frac{1}{c_4+\frac{1}{\frac{2}{c_1}}}\ge F_R(c)\\
\tilde C(c_1,c_4):=2c_1+c_4\le C(c)\end{array}$                            \\ \hline
\multicolumn{1}{l|}{$\begin{array}{l}
\mbox{\bf Case 10.2}\\
c_1>c_2+c_3\\
F=c_2+c_3+c_4\\
c_2,c_3,c_4\mbox{ are plastic}
\end{array}$              }& $\def\arraystretch{1.5}
\begin{array}{l}\tilde F_R(c_2+c_3,c_4):=0.2(c_2 +c_3+ c_4) + 0.1 \cdot\frac{1}{c_4+\frac{\frac{1}{2}}{c_2+c_3}}\ge F_R(c)\\
\tilde C(c_2+c_3,c_4):=2(c_2+c_3)+c_4\le C(c)\end{array} $                     
\end{tabular}
\caption{Summary of a one-variable estimate for $F_R$ and $C$ in the Case 10, see Appendix (Section~\ref{SectionApp}) for details.
}
\label{case10table}
\end{table}

\noindent Since formulas of Table~\ref{case7table} for Cases 10.1 and 10.2 are equivalent, it is sufficient to address just Case 10.1. 

\begin{figure}[h]
    \centering
    \vskip0.2cm
    \includegraphics[scale=0.9, trim={0cm 16cm 5cm 3cm}]{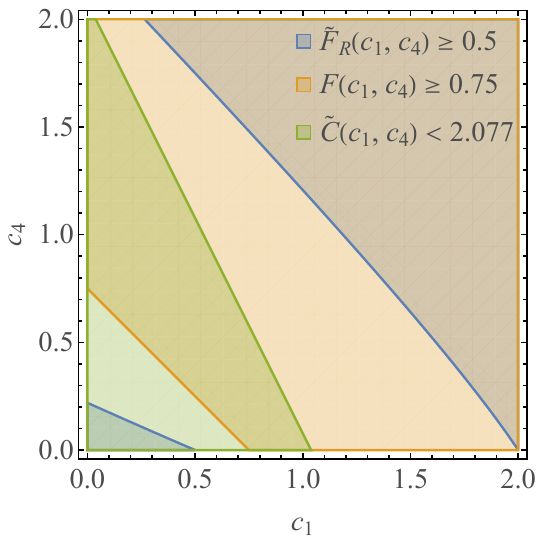} 
    
    \caption{Vizualization of the two-variable estimate for Case 10.}
    \label{case10fig}
\end{figure}
\FloatBarrier

\noindent Along the lines of the previous cases, the lack of intersection of the 3 region in Fig.~\ref{case10fig} implies that the minimal cost of fabrication in Case 10 will be higher than that of Case 9.

\section{Appendix}\label{SectionApp}

\noindent {\bf Derivation of the 2-variable estimates for $F_R$ and $C$ in the Case 9.2 (Table~\ref{case9table}).}

\vskip0.3cm

\noindent Using formula (\ref{resistance9}) for $R$ and the left column of table~\ref{case9table} for $F$, the multi-functional performance functional $F_R(c)$ computes as
$$
  F_R(c_1,c_2,c_3,c_4)=0.2(c_2 + c_4) + 0.1 \left(\frac{1}{c_1} + \frac{1}{c_4 + \frac{1}{\frac{1}{c_2} + \frac{1}{c_3}}}\right).
$$

\noindent Since $c_2<c_3$, we have 
$\frac{2}{c_2}>\frac{1}{c_2} + \frac{1}{c_3}$,
or 
$\frac{c_2}{2}<\frac{1}{\frac{1}{c_2} + \frac{1}{c_3}}$. Therefore, 
$\frac{c_2}{2}+c_4<\frac{1}{\frac{1}{c_2} + \frac{1}{c_3}}+c_4$ or 
$\frac{1}{\frac{c_2}{2}+c_4}>\frac{1}{c_4+\frac{1}{\frac{1}{c_2}+\frac{1}{c_3}}}$. On the other hand, $c_1>c_2+c_4$ implies $\frac{1}{c_2+c_4}>\frac{1}{c_1}$. Combining the last two sentences together, we get $\tilde F_R(c_2,c_4)>F_R(c_1,c_2,c_3,c_4)$ where $\tilde F_R$ is as in Table~\ref{case9table} (case 9.2).

\vskip0.2cm

\noindent To estimate the cost $C$ we will use consequently $c_2<c_3$ and $c_1>c_2+c_4$ to obtain
$$
  C(c_1,c_2,c_3,c_4)=c_1+c_2+c_3+c_4>c_1+2c_2+c_4>3c_2+2c_4, 
$$
which is the required estimate in the case 9.2 line of Table~\ref{case9table}.\\

\noindent {\bf Derivation of the 2-variable estimates for $F_R$ and $C$ in the Case 9.1 (Table~\ref{case9table}).}

\vskip0.3cm

\noindent Using formula (\ref{resistance9}) for $R$ and the left column of table~\ref{case9table} for $F$, the multi-functional performance functional $F_R(c)$ computes as
$$
  F_R(c_1,c_2,c_3,c_4)=0.2(c_1) + 0.1 \left(\frac{1}{c_1} + \frac{1}{c_4 + \frac{1}{\frac{1}{c_2} + \frac{1}{c_3}}}\right).
$$

\noindent Since $c_2<c_3$, we have 
$\frac{2}{c_2}>\frac{1}{c_2} + \frac{1}{c_3}$,
or 
$\frac{c_2}{2}<\frac{1}{\frac{1}{c_2} + \frac{1}{c_3}}$. Therefore, 
$\frac{c_2}{2}+c_4<\frac{1}{\frac{1}{c_2} + \frac{1}{c_3}}+c_4$ or 
$\frac{1}{\frac{c_2}{2}+c_4}>\frac{1}{c_4+\frac{1}{\frac{1}{c_2}+\frac{1}{c_3}}}$. On the other hand, $c_1<c_2+c_4$ implies $c_1-c_2<c_4$, or $\frac{c_2}{2}+c_1-c_2 < c_4 +\frac{c_2}{2}$, which implies $\frac{1}{c_1-\frac{c_2}{2}}>\frac{1}{c_4+\frac{c_2}{2}}$. Rewriting everything together we get, $\frac{1}{c_1-\frac{c_2}{2}}>\frac{1}{\frac{c_2}{2}+c_4}>\frac{1}{c_4+\frac{1}{\frac{1}{c_2}+\frac{1}{c_3}}}$. Adding $\frac{1}{c_1}$ on both sides we get, $\frac{1}{c_1}+\frac{1}{c_1-\frac{c_2}{2}} > \frac{1}{c_1}+\frac{1}{c_4+\frac{1}{\frac{1}{c_2}+\frac{1}{c_3}}}$, implying,  $\tilde F_R(c_2,c_2)>F_R(c_1,c_2,c_3,c_4)$ where $\tilde F_R$ is as in Table~\ref{case9table} (case 9.2).

\vskip0.2cm

\noindent To estimate the cost $C$ we will use consequently $c_2<c_3$ and $c_1<c_2+c_4$ to obtain
$$
  C(c_1,c_2,c_3,c_4)=c_1+c_2+c_3+c_4>2c_1+c_2, 
$$
which is the required estimate in the case 9.2 line of Table~\ref{case9table}.\\



\noindent \textbf{Derivation of estimates for $F_R$ and $C$ for Case 1.1 (Table \ref{case1table}).}

\vskip0.3cm

\noindent Using formula (\ref{opt1}) for Resistance and left column of Table (\ref{case1table}) for F, the multi-functional performance functional $F_R(c)$ computes as, 

$$
  F_R(c_1,c_2,c_3,c_4)=0.2c_1 + 0.1 \cdot \left (\frac{1}{c_1}+\frac{1}{c_2+c_3}+\frac{1}{c_4} \right)
$$

\noindent Since $c_1 < c_4$, we have $\frac{1}{c_1} > \frac{1}{c_4}$. Also, since $c_1 < c_2+c_3$, we have $\frac{1}{c_1} > \frac{1}{c_2+c_3}$. Combining the two together, we get $\frac{2}{c_1} > \frac{1}{c_2+c_3} + \frac{1}{c_4}$. This implies $\tilde{F}_R(c_1) > F_R (c_1,c_2,c_3,c_4)$ where $\tilde{F}_R$ is as in Table(\ref{case1table}) (Case 1.1). 

\vspace{0.2cm}

\noindent To estimate the cost $C$ we will use $c_1<c_4$ and $c_1 < c_2+c_3$ to obtain, 
$$C(c_1,c_2,c_3,c_4) = c_1+c_2+c_3+c_4 > 3c_1$$
which is the required cost estimate for Case 1.1, Table (\ref{case1table}). Similar analysis can be done for Case 1.2 following the same steps.

\vspace{0.3cm}

\vspace{0.5cm}


\noindent \textbf{Derivation of estimates for $F_R$ and $C$ for Case 2 (Table \ref{case2table}).}

\vskip0.3cm

\noindent Using formula (\ref{opt2}) for Resistance and left column of Table (\ref{case2table}) for F, the multi-functional performance functional $F_R(c)$ computes as, 

$$
  F_R(c_1,c_2,c_3,c_4)=0.2c_1 + 0.1 \cdot \left(\frac{1}{c_1}+\frac{1}{c_2}+\frac{1}{c_3}+\frac{1}{c_4}\right)
$$

\noindent Since $c_1<c_2$, $c_1<c_3$ and $c_1<c_4$, we get $\frac{1}{c_1}>\frac{1}{c_2}$, $\frac{1}{c_1}>\frac{1}{c_3}$ and $\frac{1}{c_1}>\frac{1}{c_4}$. Therefore, $\frac{4}{c_1} > \frac{1}{c_1}+\frac{1}{c_2}+\frac{1}{c_3}+\frac{1}{c_4}$ implying $\tilde{F}_R(c_1) > F_R (c_1,c_2,c_3,c_4)$.

\vspace{0.35cm}

\noindent To estimate the cost $C$ we will use $c_1<c_2$, $c_1<c_3$ and $c_1<c_4$, to obtain, 
$$C(c_1,c_2,c_3,c_4) = c_1+c_2+c_3+c_4 > 4c_1$$
which is the required cost estimate for Case 2, Table (\ref{case2table}).

\vspace{0.3cm}


\noindent \textbf{Derivation of estimates for $F_R$ and $C$ for Case 3 (Table \ref{case3table}).}

\vskip0.3cm

\noindent Using formula (\ref{opt3}) for Resistance and left column of Table (\ref{case3table}) for F, the multi-functional performance functional $F_R(c)$ computes as, 

$$
  F_R(c_1,c_2,c_3,c_4)=0.2(c_1+c_3) + 0.1 \cdot \frac{1}{\frac{1}{\frac{1}{c_1}+\frac{1}{c_2}}+\frac{1}{\frac{1}{c_3}+\frac{1}{c_4}}}
$$

\noindent Since $c_1 < c_2$, we get  $\frac{2}{c_1}>\frac{1}{c_1}+\frac{1}{c_2}$, which implies $\frac{1}{\frac{2}{c_1}}<\frac{1}{\frac{1}{c_1}+\frac{1}{c_2}}$. Similarly since $c_3<c_4$, we get, $\frac{1}{\frac{2}{c_3}} < \frac{1}{\frac{1}{c_3}+\frac{1}{c_4}}$. Combining them we get, $\frac{1}{\frac{2}{c_1}}+\frac{1}{\frac{2}{c_3}} < \frac{1}{\frac{1}{c_1}+\frac{1}{c_2}}+\frac{1}{\frac{1}{c_3}+\frac{1}{c_4}}$. Taking the reciprocal on both sides we get, $\frac{1}{\frac{c_1}{2}+\frac{c_3}{2}} > \frac{1}{\frac{1}{\frac{1}{c_1}+\frac{1}{c_2}}+\frac{1}{\frac{1}{c_3}+\frac{1}{c_4}}}$. As a result, $\tilde{F}_R(c_1,c_3)> F_R(c_1,c_2,c_3,c_4)$.

\vspace{0.35cm}

\noindent To estimate the cost $C$ we will use $c_1<c_2$ and $c_3 < c_4$ to obtain, 
$$C(c_1,c_2,c_3,c_4) = c_1+c_2+c_3+c_4 > 2c_1+2c_3$$
which is the required cost estimate for Case 3, Table (\ref{case3table}).



\noindent \textbf{Derivation of estimates for $F_R$ and $C$ for Case 4 (Table \ref{case4table}).}

\vskip0.3cm

\noindent Using formula (\ref{opt4}) for Resistance and left column of Table (\ref{case4table}) for F, the multi-functional performance functional $F_R(c)$ computes as, 

$$
  F_R(c_1,c_2,c_3,c_4)=0.2(c_1+c_4) + 0.1 \cdot \frac{1}{\frac{1}{\frac{1}{c_1}+\frac{1}{c_2}+\frac{1}{c_3}}+c_4}
$$

\noindent Since $c_1 < c_2$, and $c_2<c_3$, that is $\frac{1}{c_1} >\frac{1}{c_2}$ and $\frac{1}{c_2}>\frac{1}{c_3}$, then we have, $\frac{3}{c_1}>\frac{1}{c_1}+\frac{1}{c_2}+\frac{1}{c_3}$ which implies, $\frac{1}{\frac{3}{c_1}}+c_4 < \frac{1}{\frac{1}{c_1}+\frac{1}{c_2}+\frac{1}{c_3}}+c_4$. Taking the reciprocal once again we get, 
$\frac{1}{\frac{1}{\frac{3}{c_1}}+c_4}>\frac{1}{\frac{1}{\frac{1}{c_1}+\frac{1}{c_2}+\frac{1}{c_3}}+c_4}$. Thus, $\tilde{F}_R(c_1,c_4)> F_R (c_1,c_2,c_3,c_4)$.

\vspace{0.35cm}

\noindent To estimate the cost $C$ we will use $c_1<c_2$ and $c_1 < c_3$ to obtain, 
$$C(c_1,c_2,c_3,c_4) = c_1+c_2+c_3+c_4 > 3c_1+c_4$$
which is the required cost estimate for Case 4, Table (\ref{case4table}).


\vspace{0.3cm}


\noindent \textbf{Derivation of estimates for $F_R$ and $C$ for Case 5 (Table \ref{case5table}).}

\vskip0.3cm

\noindent Using formula (\ref{opt5}) for Resistance and left column of Table (\ref{case5table}) for F, the multi-functional performance functional $F_R(c)$ computes as, 

$$
  F_R(c_1,c_2,c_3,c_4)=0.2(c_1+c_3+c_4) + 0.1 \cdot \frac{1}{\frac{1}{\frac{1}{c_1}+\frac{1}{c_2}}+c_3+c_4}
$$

\noindent Since $c_1 < c_2$, we get $\frac{1}{\frac{2}{c_1}} < \frac{1}{\frac{1}{c_1}+\frac{1}{c_2}}$, which implies, $\frac{1}{\frac{2}{c_1}} +c_3+c_4 < \frac{1}{\frac{1}{c_1}+\frac{1}{c_2}}+c_3+c_4$. Taking the reciprocal of the term we get, $\frac{1}{\frac{1}{\frac{2}{c_1}} +c_3+c_4} > \frac{1}{\frac{1}{\frac{1}{c_1}+\frac{1}{c_2}}+c_3+c_4}$. This concludes that, $\tilde{F}_R(c_1,c_3,c_4) > F_R(c_1,c_2,c_3,c_4)$.

\vspace{0.35cm}

\noindent To estimate the cost $C$ we will use $c_1<c_2$ to obtain, 
$$C(c_1,c_2,c_3,c_4) = c_1+c_2+c_3+c_4 > 2c_1+c_3+c_4$$
which is the required cost estimate for Case 5, Table (\ref{case5table}).

\vspace{0.3cm}


\noindent \textbf{Derivation of estimates for $F_R$ and $C$ for Case 6.1 (Table \ref{case6table}).}

\vskip0.3cm

\noindent Using formula (\ref{opt6}) for Resistance and left column of Table (\ref{case6table}) for F, the multi-functional performance functional $F_R(c)$ computes as, 

$$
  F_R(c_1,c_2,c_3,c_4)=0.2c_1 + 0.1 \cdot \left(\frac{1}{c_1}+\frac{1}{c_2+c_3+c_4}\right)
$$

\noindent Since $c_1 < c_2+c_3+c_4$, we get $\frac{1}{c_1}>\frac{1}{c_2+c_3+c_4}$, so, $\frac{2}{c_1}>\frac{1}{c_1}+\frac{1}{c_2+c_3+c_4}$. Therefore, $\tilde{F}_{R}(c_1) > F_R(c_1,c_2,c_3,c_4)$. 

\vspace{0.35cm}

\noindent To estimate the cost $C$ we will use $c_1<c_2+c_3+c_4$ to obtain, 
$$C(c_1,c_2,c_3,c_4) = c_1+c_2+c_3+c_4 > 2c_1$$
which is the required cost estimate for Case 6.1, Table (\ref{case6table}).Similar analysis can be done
for Case 6.2 following the same steps. 


\vspace{0.3cm}


\noindent \textbf{Derivation of estimates for $F_R$ and $C$ for Case 7.1 (Table \ref{case7table}).}

\vskip0.3cm

\noindent Using formula (\ref{opt7}) for Resistance and left column of Table (\ref{case7table}) for F, the multi-functional performance functional $F_R(c)$ computes as, 

$$
  F_R(c_1,c_2,c_3,c_4)=0.2(c_1+c_3) + 0.1 \cdot \left(\frac{1}{c_1+c_3}+\frac{1}{c_2+c_4}\right)
$$

\noindent Since $c_1+c_3 < c_2+c_4$, we get $\frac{1}{c_1+c_3} > \frac{1}{c_2+c_4}$, which implies, $\frac{2}{c_1+c_3}> \frac{1}{c_1+c_3}+\frac{1}{c_2+c_4}$. Therefore, $\tilde{F}_R(c_1,c_3) = F_R(c_1,c_2,c_3,c_4)$.
\vspace{0.35cm}

\noindent To estimate the cost $C$ we will use $c_1+c_3<c_2+c_4$ to obtain, 
$$C(c_1,c_2,c_3,c_4) = c_1+c_2+c_3+c_4 > 2(c_1+c_3)$$
which is the required cost estimate for Case 7.1, Table (\ref{case7table}).Similar analysis can be done for Case 7.2 following the same steps. Similar analysis can be done
for Case 7.2 following the same steps.

\vspace{0.3cm}


\noindent \textbf{Derivation of estimates for $F_R$ and $C$ for Case 8 (Table \ref{case8table}).}

\vskip0.3cm

\noindent Using formula (\ref{opt8}) for Resistance and left column of Table (\ref{case8table}) for F, the multi-functional performance functional $F_R(c)$ computes as, 

$$
  F_R(c_1,c_2,c_3,c_4)=0.2(c_1+c_2+c_3+c_4) + 0.1 \cdot \left(\frac{1}{c_1+c_2+c_3+c_4}\right)
$$



$$C(c_1,c_2,c_3,c_4) = c_1+c_2+c_3+c_4 $$

\vspace{0.3cm}


\noindent \textbf{Derivation of estimates for $F_R$ and $C$ for Case 10.1} (Table \ref{case10table}).

\vskip0.3cm

\noindent Using formula (\ref{opt10}) for Resistance and left column of Table (\ref{case10table}) for F, the multi-functional performance functional $F_R(c)$ computes as, 

$$
  F_R(c_1,c_2,c_3,c_4)=0.2(c_1+c_4) + 0.1 \cdot \left(\frac{1}{c_4+\frac{1}{\frac{1}{c_1}+\frac{1}{c_2+c_3}}}\right)
$$

\noindent Since $c_1< c_2+c_3$, we get $\frac{1}{c_1} > \frac{1}{c_2+c_3}$, which implies, $\frac{2}{c_1}> \frac{1}{c_1}+\frac{1}{c_2+c_3}$. From this, we can obtain, $c_4 + \frac{1}{\frac{2}{c_1}} < c_4 +\frac{1}{\frac{1}{c_1}+\frac{1}{c_2+c_3}}$. Taking the reciprocal we get, $\frac{1}{c_4 + \frac{1}{\frac{2}{c_1}}}>\frac{1}{c_4 +\frac{1}{\frac{1}{c_1}+\frac{1}{c_2+c_3}}}$. Therefore, $\tilde{F}_R(c_1,c_4)=F_R(c_1,c_2,c_3,c_4)$.
\vspace{0.35cm}

\noindent To estimate the cost $C$ we will use $c_1<c_2+c_3$ to obtain, 
$$C(c_1,c_2,c_3,c_4) = c_1+c_2+c_3+c_4 > 2c_1+c_4$$
which is the required cost estimate for Case 10.1, Table (\ref{case10table}). Similar analysis can be done for Case 10.2 following the same steps. 

\section*{Data availability statement}

\noindent All data that support the findings of this study are included within the article (and any
supplementary files).

\end{document}